\documentclass[preprint,12pt]{elsarticle}

\usepackage{amssymb}
\usepackage{amsmath}

\begin{document}

\makeatletter
\def\ps@pprintTitle{%
 \let\@oddhead\@empty
 \let\@evenhead\@empty
 \let\@oddfoot\@empty
 \let\@evenfoot\@empty}
\makeatother

\begin{frontmatter}



\begin{center}
    {\LARGE \textbf{Preprint Version}} \\[1em]
    \textit{This is a preprint version.} \\[1em]
    Date: \today
\end{center}

\title{Capacity Constraints in Ball and Urn Distribution Problems}


\author[Label1]{Jingwei Li}
\ead{jingweili2017@gmail.com}
\author[Label1]{Thomas G. Robertazzi\corref{cor1}}
\ead{thomas.robertazzi@stonybrook.edu}
\cortext[cor1]{Corresponding author}

\affiliation[Label1]{organization={Department of Electrical and Computer Engineering,
Stony Brook University},
            addressline={100 Nicolls Rd},
            city={Stony Brook},
            postcode={11790},
            state={NY},
            country={USA}}
            
\date{\textbf{Preprint version – \today}}

\begin{abstract}
This paper explores the distribution of indistinguishable balls into distinct urns with varying capacity constraints, a foundational issue in combinatorial mathematics with applications across various disciplines.  We present a comprehensive theoretical framework that addresses both upper and lower capacity constraints under different distribution conditions, elaborating on the combinatorial implications of such variations. Through rigorous analysis, we derive analytical solutions that cater to different constrained environments, providing a robust theoretical basis for future empirical and theoretical investigations. These solutions are pivotal for advancing research in fields that rely on precise distribution strategies, such as physics and parallel processing. The paper not only generalizes classical distribution problems but also introduces novel methodologies for tackling capacity variations, thereby broadening the utility and applicability of distribution theory in practical and theoretical contexts.
\end{abstract}


\begin{highlights}
\item Developed a comprehensive framework for distributing indistinguishable balls into urns with varying capacity constraints, enhancing the classical understanding of combinatorial distribution problems.

\item Derived novel analytical solutions for different constrained environments, providing essential tools for future research in fields requiring precise distribution strategies, such as physics and parallel processing.

\item Introduced innovative methodologies for analyzing upper and lower capacity constraints in urn models, significantly broadening the applicability of distribution theory to complex real-world scenarios.
\end{highlights}


\end{frontmatter}



\section{Introduction}
\label{sec1}

In the context of combinatorial mathematics, the distribution problem is a significant and foundational topic. Understanding how to allocate items into various spaces has numerous meaningful applications across diverse research fields. For instance, in the partial derivatives problem, given an analytical function $f(x_1, x_2, \dots, x_n)$, the number of distinct partial derivatives of order $r$ is $\binom{n+r-1}{r}$ \cite{feller1991introduction}. 

In physics, the Maxwell-Boltzmann distribution can be derived using the "Boltzmann relation" between entropy $S$ and the probability of occupation $W$. Specifically, $W$ is derived as $W = \frac{1}{\prod_k N_k!},$ 
using the combinatorial distribution method when $N$ indistinguishable particles are distributed into $k$ occupation cells, where $N_k$ represents the number of particles in each cell \cite{gordon2002maxwell}. Furthermore, the combinatorial distribution method is instrumental in calculating the distribution of bosons, as for bosons, the number of particles in a single cell is constrained only by the total number of particles \cite{gazdzicki2020notes}. 

In the domain of parallel processing, particularly in divisible load scheduling \cite{ko2008signature, kyong2012greedy}, the combinatorial distribution method can be employed to calculate the expected signature search time at each layer of multi-level tree networks when the total number of signatures is known per layer \cite{ying2013signature}.

Given the importance of the combinatorial distribution problem, there is a need to extend its study to scenarios where urns (or boxes or cells) have capacity constraints under varying conditions. Additionally, the relationship between the number of balls and urns introduces new challenges and variations in the distribution process.

This research investigates the problem of distributing indistinguishable (or identical) balls into distinct boxes while considering different capacity constraints. A systematic analysis of the problem is presented. In Section 2, we provide a detailed examination of multiple constraints on box capacity under the condition where the number of balls exceeds the number of boxes. In Section 3, we analyze scenarios where the number of balls is smaller than the number of boxes, with varying capacity constraints. The capacity of each box is divided into distinct segments that form a complete capacity-constrained space. The conclusions of this study are presented in Section 4.




\section{Derivations of Different Scenarios When $m \geq n$}
\label{sec2}

Based on the "Stars and Bars" method \cite{feller1991introduction}, when $m \geq n$, there are $\binom{m-1}{n-1}$ ways to distribute $m$ identical balls into $n$ distinct boxes if no box is empty. However, if empty boxes are allowed, there are $\binom{m+n-1}{n-1}$ possible ways to distribute the balls. When $m < n$, some boxes will inevitably remain empty, and the total number of combinations is $\binom{m+n-1}{n-1}$.

In the following sections, we discuss distribution problems under additional constraints on box capacities. Specifically, we assume all balls must be distributed with no remainder, and each box has the same capacity $k$, where $k_1 \leq k \leq k_2$. Analytical solutions are derived for various ranges of $k_1$ and $k_2$.

\subsection{$k_2 \geq m$}
\label{sec21}
When $k_2 \geq m$, the maximum capacity constraint is effectively irrelevant, as no box will hold more than $m$ balls. Thus, the problem simplifies significantly.

\subsubsection{$k_1 = 1$ and $k_2 \geq m$}
\label{sec211}
This scenario is equivalent to distributing $m$ identical balls into $n$ distinct boxes with no empty boxes allowed. The total number of possible distributions is $\binom{m-1}{n-1}$.

\subsubsection{$k_1 = 0$ and $k_2 \geq m$}
\label{sec212}
In this case, empty boxes are allowed, making the scenario identical to distributing $m$ identical balls into $n$ distinct boxes without restrictions on empty boxes. The total number of possible distributions is $\binom{m+n-1}{n-1}$.

\subsubsection{$ k_1 > \lfloor \frac{m}{n} \rfloor $ and $k_2 \geq m$} \label{213}
The number of combinations is zero since the constraint $k_1 > \lfloor \frac{m}{n} \rfloor$ is invalid. 

\textbf{Lemma 2.1.} \textit{When distributing $m$ identical balls into $n$ distinct boxes, if the box capacity $k$ has the constraint $k \geq k_1$, then $k_1$ cannot exceed $\lfloor \frac{m}{n} \rfloor$.}

\textbf{Proof.} There are $m$ balls to distribute among $n$ boxes, with $\lfloor \frac{m}{n} \rfloor \leq \frac{m}{n} < \lfloor \frac{m}{n} \rfloor + 1$. If $k_1 > \lfloor \frac{m}{n} \rfloor$, then $k > \lfloor \frac{m}{n} \rfloor$ because $k \geq k_1$. Assuming $k = \lfloor \frac{m}{n} \rfloor + 1$, each box would require exactly $\lfloor \frac{m}{n} \rfloor + 1$ balls. The total number of balls required would then be $n \cdot (\lfloor \frac{m}{n} \rfloor + 1)$, which exceeds the total $m$ balls available, violating the initial condition. Therefore, the constraint $k_1 > \lfloor \frac{m}{n} \rfloor$ is invalid, and the number of possible distributions is zero.

\subsubsection{$ 1 \leq k_1 \leq \lfloor \frac{m}{n} \rfloor $ and $k_2 \geq m$} \label{214}
Based on the discussion in Section \ref{213}, the constraint $1 \leq k_1 \leq \lfloor \frac{m}{n} \rfloor$ is valid. The problem is to distribute $m$ identical balls into $n$ distinct boxes such that each box contains at least $k_1$ balls. 

To solve this, a two-stage method is employed. In stage 1, $k_1$ balls are placed in each box. This consumes $n \cdot k_1$ balls, leaving $m - n \cdot k_1$ balls. In stage 2, the remaining $m - n \cdot k_1$ balls are distributed freely among the $n$ boxes, where each box can receive zero or more balls. 

There is only one way to perform stage 1 since the balls are identical. Based on the discussion in Section \ref{sec212}, the number of ways to complete stage 2 is $\binom{m - n \cdot k_1 + n - 1}{n - 1}$. Thus, the total number of ways to distribute the balls is $1\cdot \binom{m - n \cdot k_1 + n - 1}{n - 1} = \binom{m - n \cdot k_1 + n - 1}{n - 1}$.

\subsection{$k_1=0$}\label{sec22}
When $k_1=0$, the lower bound on box capacity is removed, meaning the number of balls in each box can only be constrained by the upper bound $k_2$. This simplifies the analysis to focus solely on the upper capacity constraint.

\subsubsection{$k_1 = 0$ and $ \lfloor \frac{m}{2} \rfloor \leq k_2 < m$} \label{221}
Conditioning on  $k_2 < m$, the complement of the constraint $k \leq k_2$ is $k > k_2$, equivalent to $ \lfloor \frac{m}{2} \rfloor + 1\leq k_2+1 \leq k \leq m - 1$. Since there are only $m$ balls to distribute, at most one box can contain more than $\lfloor \frac{m}{2} \rfloor$ balls. Consequently, the event where at least one box contains more than $\lfloor \frac{m}{2} \rfloor$ balls is equivalent to the event where exactly one box exceeds this limit.

To solve for the number of ways to distribute $m$ balls among $n$ boxes under this condition, a two-stage method is employed:
\begin{itemize}
    \item \textbf{Stage 1}: Select one of the $n$ boxes to hold the $k_2 + 1$ balls. The number of combinations for this step is  $\binom{n}{1}$.
    \item \textbf{Stage 2}: Distribute the remaining $m - (k_2 + 1)$ balls among the $n$ boxes without capacity constraints. Based on the result of Section \ref{sec212}, the number of combinations for this step is $\binom{m - k_2 - 1 + n - 1}{n-1}$.
\end{itemize}

Thus, the total number of combinations for this scenario is:
\[  \binom{n}{1} \cdot \binom{m - k_2 - 1 + n - 1}{n-1} . \].

Returning to the original problem, the total sample space corresponds to distributing $m$ balls into $n$ boxes without a lower limit on the capacity of the box. The total number of combinations in this sample space is $\binom{m+n-1}{n-1}$. Therefore, the number of ways to distribute $m$ balls under the conditions $k_1 = 0$ and $ \lfloor \frac{m}{2} \rfloor \leq k_2 < m$ is given by:
\[
\binom{m + n - 1}{n-1} - \binom{n}{1} \cdot \binom{m - k_2 + n - 2}{n-1}.
\]

\subsubsection{$k_1 = 0$ and $k_2 < \lceil \frac{m}{n} \rceil$} \label{222}

It is an impossible constraint if $k_2 < \lceil \frac{m}{n} \rceil$, as the condition renders the distribution of $m$ balls into $n$ boxes infeasible. 

\textbf{Lemma 2.2.} \textit{When distributing $m$ identical balls into $n$ distinct boxes, if the box capacity $k$ is constrained such that $k \leq k_2$, the value of $k_2$ cannot be smaller than $\lceil \frac{m}{n} \rceil$.}

\textbf{Proof.} If $k_2 < \lceil \frac{m}{n} \rceil$, then it follows that $n k_2 < n \lceil \frac{m}{n} \rceil$. Expanding this inequality, $n k_2 \leq n (\lceil \frac{m}{n} \rceil - 1)$ also holds. Further, we observe that:
\[
n (\lceil \frac{m}{n} \rceil - 1) < n (\frac{m}{n} +1 - 1) = n \cdot \frac{m}{n} = m.
\]
Thus, $n k_2 < m$. 

This implies that even if each box is filled to its maximum capacity, i.e., $k_2$, the total number of balls distributed is strictly less than $m$. Consequently, there would still be $m - n k_2$ remaining balls, which is not permissible since all balls must be distributed among the boxes. While empty boxes are allowed, remaining balls are not permitted under the constraints of this problem.

Therefore, the capacity $k_2$ must satisfy the condition $k_2 \geq \lceil \frac{m}{n} \rceil$. Any value of $k_2$ smaller than this threshold violates the feasibility of the distribution. 

As a result, the total number of ways to distribute the balls under the constraint $k_2 < \lceil \frac{m}{n} \rceil$ is 0.

\subsubsection{$  k_1 = 0$ and $ \lceil \frac{m}{n} \rceil \leq k_2 < \lfloor \frac{m}{2} \rfloor$} \label{223}
Assume that $X_1, X_2, \ldots, X_n$ represent the number of balls in each box, numbered from $1$ to $n$. The complement of the condition $X_i \leq \kappa$ is $X_i \geq \kappa + 1$, where $\kappa$ is a positive integer. Let event $A$ denote the distribution of $m$ identical balls into $n$ distinct boxes, such that each box contains at most $\kappa$ balls. The complementary event $A^C$ then represents the scenario where at least one box contains more than $\kappa$ balls.

Unlike Section~\ref{221}, the expression $\binom{n}{1} \cdot \binom{m - \kappa - 1 + n - 1}{n-1}$
cannot be used here to compute the number of combinations, as there may be multiple boxes containing more than $\kappa$ balls when $\lceil \frac{m}{n} \rceil \leq \kappa < \lfloor \frac{m}{2} \rfloor$. Directly applying the "Stars and Bars" method leads to repeated calculations, which would result in incorrect outcomes.

To simplify the analysis, we can impose the capacity constraint on a single specific box at a time. For example, calculating the combinations when only box $1$ has at most $\kappa$ balls is straightforward. Motivated by this idea, We propose an approach based on the Inclusion–Exclusion Principle \cite{ryser1963combinatorial} to solve this problem.


Consider event $A_1$ as distributing $m$ identical balls into $n$ distinct boxes such that $X_1 \leq \kappa$. Its complementary event, $A_1^C$, represents the scenario where $X_1 \geq \kappa + 1$. Similarly, event $A_i$ corresponds to $X_i \leq \kappa$, with its complement $A_i^C$ defined as $X_i \geq \kappa + 1$ for $i = 1, 2, \ldots, n$. The number of combinations for $A_i^C$ is $\binom{m - \kappa - 1 + n - 1}{n-1},$
and for $A_i$ it is 
$\binom{m + n - 1}{n-1} - \binom{m - \kappa - 1 + n - 1}{n-1}.$

In set theory and Boolean algebra \cite{goodstein2007boolean}, the principle of "union and intersection interchange under complementation" is often stated as follows:
\[
(A \cap B)^C = A^C \cup B^C,
\]
where $A$ and $B$ are two events. This formula can be generalized for multiple events as:
\[
\left( \bigcap_{i=1}^n A_i \right)^C = \bigcup_{i=1}^n A_i^C.
\]
Equivalently, we have:
\[
\bigcap_{i=1}^n A_i = \left( \bigcup_{i=1}^n A_i^C \right)^C.
\]
This is commonly known as "De Morgan's Laws" \cite{goodstein2007boolean}. Using this framework, we proceed to compute the probability $P\left(\bigcup_{i=1}^n A_i^C\right)$.

For $n = 2$, the probability of the union of two events is:
\begin{equation}
    P(A_1^C \cup A_2^C) = P(A_1^C) + P(A_2^C) - P(A_1^C \cap A_2^C).
\end{equation}

For $n$ events, the general formula is given by:
\begin{equation} \label{union}
\begin{split}
    P\left(\bigcup_{i=1}^n A_i^C\right) = & \sum_{i=1}^n P(A_i^C) - \sum_{i_1 \neq i_2} P(A_{i_1}^C \cap A_{i_2}^C) \\
    & + \sum_{i_1 \neq i_2 \neq i_3} P(A_{i_1}^C \cap A_{i_2}^C \cap A_{i_3}^C) - \cdots \\
    & + (-1)^{n+1}P(A_{i_1}^C \cap A_{i_2}^C\cap 
 A_{i_3}^C...\cap A_{i_{n-1}}^C\cap A_{i_{n}}^C).
\end{split}
\end{equation}

Since $P(A_1^C) = P(A_2^C) = \cdots = P(A_n^C)$, Equation~\eqref{union} simplifies to:
\begin{equation} \label{prob_sum}
\begin{split}
    P\left(\bigcup_{i=1}^n A_i^C\right) = & \binom{n}{1} P(A_i^C) - \binom{n}{2} P(A_1^C \cap A_2^C) \\
    & + \binom{n}{3} P(A_1^C \cap A_2^C \cap A_3^C) - \cdots \\
    & + (-1)^{r+1}\binom{n}{r}P(A_{1}^C \cap A_{2}^C\cap A_{3}^C...\cap A_{r}^C) + ... \\
    & + (-1)^{n+1}P(A_{1}^C \cap A_{2}^C\cap A_{3}^C...A_{{n}}^C)
\end{split}
\end{equation}

Notably, the summation in Equation~\eqref{prob_sum} terminates at $\lfloor \frac{m}{\kappa+1} \rfloor$, as it is the maximum number of boxes that can contain more than $\kappa$ balls.


\textbf{Lemma 2.3. } \label{lemma23} \textit{ When distributing $m$ identical balls into $n$ distinct boxes, the maximum number of boxes that can contain more than $\kappa$ balls is $\lfloor \frac{m}{\kappa+1} \rfloor$.}

\textbf{Proof.} Suppose there are more than $\lfloor \frac{m}{\kappa+1} \rfloor$ boxes with more than $\kappa$ balls. Let this number be $\lfloor \frac{m}{\kappa+1} \rfloor + 1$. Assume that each of these boxes contains exactly $\kappa+1$ balls. The total number of balls required would then be:
\[
(\kappa+1) \cdot \left(\lfloor \frac{m}{\kappa+1} \rfloor + 1\right) > (\kappa+1) \cdot \frac{m}{\kappa+1} = m.
\]
This contradicts the assumption that only $m$ balls are available for distribution. Hence, the maximum number of boxes with more than $\kappa$ balls is $\lfloor \frac{m}{\kappa+1} \rfloor$. 



Using \textbf{Lemma 2.3}, the probability of $\bigcup_{i=1}^n A_i^C$ becomes:

\begin{equation} 
    \begin{split}
        P\left(\bigcup_{i=1}^{n} A_i^C\right)
        =& \binom{n}{1} P(A_i^C) - \binom{n}{2} P(A_1^C \cap A_2^C) + \binom{n}{3} P(A_{1}^C \cap A_{2}^C \cap A_{3}^C) -  \\
        \dots &+ (-1)^{r+1} \binom{n}{r} P(A_{1}^C \cap A_{2}^C \cap A_{3}^C \dots \cap A_{r}^C) + \\
        \dots &+ (-1)^{\lfloor \frac{m}{\kappa+1} \rfloor+1} \binom{n}{\lfloor \frac{m}{\kappa+1} \rfloor} P(A_{1}^C \cap A_{2}^C \cap A_{3}^C \dots \cap A_{{\lfloor \frac{m}{\kappa+1} \rfloor}}^C).
    \end{split}
\end{equation}


Further, the probabilities for all the different intersected events from $A_1^C$ to $A_n^C$ are expressed as follows:


\begin{equation}
    \begin{cases}
      P(A_i^C)= 
      
      \frac{\binom{m-(\kappa+1)+n-1}{n-1}}{\binom{m+n-1}{n-1}}, &i=1,2,...,n\\ 
      P(A_i^C \cap A_j^C)=\frac{\binom{m-2(\kappa+1)+n-1}{n-1}}{\binom{m+n-1}{n-1}},&i\ne j, i\&j=1,2,...,n\\
      P(A_{i_1}^C ...\cap A_{i_j}^C...\cap A_{i_r}^C)=\frac{\binom{m-r(\kappa+1)+n-1}{n-1}}{\binom{m+n-1}{n-1}},
      &i_1\ne i_2 \ne...i_r, r=3,4,...,\\
      &{\lfloor \frac{m}{\kappa+1} \rfloor},i_j=1,2,...,n\\
    \end{cases}  
\end{equation}

Thus, the probability of $\bigcup_{i=1}^n A_i^C$ can be written as:
\begin{equation}
    P\left(\bigcup_{i=1}^n A_i^C\right) = \frac{\sum_{\alpha=1}^{\lfloor \frac{m}{\kappa+1} \rfloor} (-1)^{\alpha+1} \binom{n}{\alpha} \binom{m-\alpha \cdot (\kappa+1)+n-1}{n-1}}{\binom{m+n-1}{n-1}}.
\end{equation}

By De Morgan's Law, the probability of $\bigcap_{i=1}^n A_i$ is given by:
\begin{equation}
\begin{split}
     P\left(\bigcap_{i=1}^n A_i\right) &= P\left(\left(\bigcup_{i=1}^n A_i^C\right)^C\right) \\
     &= 1 - P\left(\bigcup_{i=1}^n A_i^C\right) \\
     &= 1 - \frac{\sum_{\alpha=1}^{\lfloor \frac{m}{\kappa+1} \rfloor} (-1)^{\alpha+1} \binom{n}{\alpha} \binom{m-\alpha \cdot (\kappa+1)+n-1}{n-1}}{\binom{m+n-1}{n-1}} \\
     &= P_m^{n, \kappa}.
\end{split}
\end{equation}

We define a new symbol $P_m^{n, \kappa}$ to represent the probability of event $E$, where $E$ is the event of distributing $m$ identical balls into $n$ distinct boxes such that each box contains at most $\kappa$ balls, with the constraint $\lceil \frac{m}{n} \rceil \leq \kappa < m$. Furthermore, we introduce $\Omega_m^{n, \kappa}$ to denote the total number of valid distributions corresponding to event $E$, which is expressed as:

\begin{equation}
    \Omega_m^{n, \kappa} = \binom{m+n-1}{n-1} - \sum_{\alpha=1}^{\lfloor \frac{m}{\kappa+1} \rfloor} (-1)^{\alpha+1} \binom{n}{\alpha} \binom{m-\alpha \cdot (\kappa+1)+n-1}{n-1}.
\end{equation}

Therefore, for the case where $k_1 = 0$ and $\lceil \frac{m}{n} \rceil \leq k_2 < \lfloor \frac{m}{2} \rfloor$, the total number of combinations to distribute $m$ identical balls into $n$ distinct boxes is:
\begin{equation}
    \Omega_m^{n, k_2} = \binom{m+n-1}{n-1} - \sum_{\alpha=1}^{\lfloor \frac{m}{k_2+1} \rfloor} (-1)^{\alpha+1} \binom{n}{\alpha} \binom{m-\alpha \cdot (k_2+1)+n-1}{n-1}.
\end{equation}

Further, if $\lfloor \frac{m}{2} \rfloor \leq k_2 < m$, the value of $\Omega_m^{n, k_2}$ reduces to:
$\binom{m + n - 1}{n - 1} - \binom{n}{1} \cdot \binom{m - k_2 + n - 2}{n - 1}$,
which corresponds to the case discussed in Section~\ref{221}.

\subsection{$1\leq k_1 \leq \lfloor \frac{m}{n} \rfloor$ and $ \lceil \frac{m}{n} \rceil \leq k_2 < m$}\label{sec3} In sections \ref{sec21} and \ref{sec22}, we discussed scenarios where the variable $k$ is constrained by only one boundary—either a lower or an upper limit, respectively. In this section, the details will be discussed when $k$ has the lower bound and upper bound at the same time. When $k_1 \ne 0$, the valid range of $k_1$ should be $1\leq k_1 \leq \lfloor \frac{m}{n} \rfloor$; when $k_2\not \geq m $, the valid range of $k_2$ should be $\lceil \frac{m}{n} \rceil \leq k_2 < m$.

Similar to the discussion in section \ref{214}, a two-stage method will be used here. In stage 1, each of the $n$ box will be distributed with $k_1$ balls. After stage 1, there are $m^* = m-n\cdot k_1$ balls remaining. In stage 2, the problem  is equal to distribute $m^*$ identical balls into $n$ distinct boxes where the box capacity $k^*$ has the constraint of $k_1^* = 0 \leq k^* \leq k_2^*$, where $k_2^*$ is a whole number. Based on the discussion in section \ref{sec2}, the valid range of $k_2^*$ should be $\lfloor \frac{m^*}{2} \rfloor \leq k_2^* < m^{*}$ or $ \lceil \frac{m^*}{n} \rceil \leq k_2^* <\lfloor \frac{m^*}{2} \rfloor$.


\subsubsection{$1 \leq k_1 \leq \lfloor \frac{m}{n} \rfloor$ and $\lfloor \frac{m}{2} \rfloor \leq k_2 < m$}\label{sec231}

If \( \lfloor \frac{m}{2} \rfloor \leq k_2 < m \), then \( \lfloor \frac{m}{2} \rfloor - k_1 \leq k_2^* = k_2 - k_1 < m - k_1 \). Simplifying further, we obtain $\lfloor \frac{m - 2k_1}{2} \rfloor \leq k_2^* < m - k_1.$

The left-hand side implies $k_2^* \geq \lfloor \frac{m}{2} \rfloor - k_1 = \lfloor \frac{m - 2k_1}{2} \rfloor \geq \lfloor \frac{m - n \cdot k_1}{2} \rfloor = \lfloor \frac{m^*}{2} \rfloor$
for \( n > 1 \); while the right-hand side implies $k_2^* < m^* = m - n \cdot k_1$ or $m^* \leq k_2^* < m - k_1$, 
since \( m - n \cdot k_1 \leq m - k_1 \).

Thus, the range of \( k_2^* \) can be rewritten as:
\[
\quad m^* \leq k_2^* < m - k_1 \quad \text{or} \quad
\lfloor \frac{m}{2} \rfloor-k_1 \leq k_2^* < m^*.
\]

\textbf{Case A: $k_1^{*} = 0$ and $m^* \leq k_2^* < m - k_1$}.  
From the analysis in Section \ref{sec212}, the number of combinations is:
\[
\binom{m^* + n - 1}{n-1} = \binom{m - nk_1 + n - 1}{n-1}.
\]

\textbf{Case B: $k_1^{*} = 0$ and $\lfloor \frac{m}{2} \rfloor - k_1 \leq k_2^{*} < m^*$}.  
The range $\lfloor \frac{m}{2} \rfloor - k_1 \leq k_2^* < m^*$ is a subset of the range $\lfloor \frac{m^*}{2} \rfloor \leq k_2^* < m^*$. Based on the discussion in \ref{221}, the number of combinations in this range is:
\[
\binom{m^* + n - 1}{n-1}-\binom{n}{1} \cdot \binom{m^* - k_2^* + n - 2}{n-1}
\]
Rewriting using $k_2^* = k_2 - k_1$ and $m^*=m-nk_1$, we get:
\[
\binom{m - n k_1 + n - 1}{n-1} - n\cdot \binom{m - n k_1 - (k_2 - k_1) + n - 2}{n-1}
\]

\subsubsection{$1 \leq k_1 \leq \lfloor \frac{m}{n} \rfloor$ and $\lceil \frac{m}{n} \rceil \leq k_2 < \lfloor \frac{m}{2} \rfloor$}\label{sec232}

For $k_2^* = k_2 - k_1$, the upper bound satisfies $k_2^* < \lfloor \frac{m}{2} \rfloor - k_1$. From Section \ref{sec231}, $\lfloor \frac{m^*}{2} \rfloor \leq \lfloor \frac{m}{2} \rfloor - k_1$. Thus, $k_2^* < \lfloor \frac{m}{2} \rfloor - k_1$ is equivalent to:
\[
\{k_2^* < \lfloor \frac{m^*}{2} \rfloor\} \cup \{\lfloor \frac{m^*}{2} \rfloor \leq k_2^* < \lfloor \frac{m}{2} \rfloor - k_1\}.
\]
For the lower bound, $k_2 \geq \lceil \frac{m}{n} \rceil$ implies $k_2^* = k_2 - k_1 \geq \lceil \frac{m}{n} \rceil - k_1=\lceil \frac{m-nk_1}{n} \rceil$, which is equivalent to $k_2^* \geq \lceil \frac{m^*}{n} \rceil$.

Thus, two cases for $k_2^*$ arise:  

\textbf{Case A: $k_1^*=0 \ and \ \lfloor \frac{m^*}{2} \rfloor \leq k_2^* < \lfloor \frac{m}{2} \rfloor - k_1$}.  
The upper bound of $k_2^*$ can either greater or smaller than $m^*$ dependent on the specific values of $m,n,k_1$ and $k_2$. When $\lfloor \frac{m}{2} \rfloor - k_1 < m^*$,  
the range  $\lfloor \frac{m^*}{2} \rfloor \leq k_2^* < \lfloor \frac{m}{2} \rfloor - k_1$ is a subset of the range $\lfloor \frac{m^*}{2} \rfloor \leq k_2^* < m^*$. From section \ref{sec231}, the number of combinations is 
\[
\binom{m - n k_1 + n - 1}{n-1}- n \cdot \binom{m - n k_1 - (k_2 - k_1) + n - 2}{n-1}
\]

When $\lfloor \frac{m}{2} \rfloor - k_1 >= m^*$, the range of $k_2^*$ turns to $\lfloor \frac{m^*}{2} \rfloor \leq k_2^* < m^*$ or $ m^* \leq k_2^* \leq  \lfloor \frac{m}{2} \rfloor - k_1$. For the former, the number of combinations is
\[
\binom{m - n k_1 + n - 1}{n-1} - n \cdot \binom{m - n k_1 - (k_2 - k_1) + n - 2}{n-1};
\]
for the latter, based on the discussion in \ref{sec212}, the number of combinations is 
\[
\binom{m^*+n-1}{n-1} = \binom{m-nk_1+n-1}{n-1}
\]

\textbf{Case B: $ k_1^*=0 \ and \ \lceil \frac{m^*}{n} \rceil \leq k_2^* < \lfloor \frac{m^*}{2} \rfloor$}.  
From Section \ref{223}, the total number of combinations is:
\[
\Omega_{m^*}^{n, k_2^*} = \binom{m^* + n - 1}{n - 1} - \sum_{\alpha = 1}^{\lfloor \frac{m^*}{k_2^* + 1} \rfloor} (-1)^{\alpha + 1} \binom{n}{\alpha} \binom{m^* - \alpha \cdot (k_2^* + 1) + n - 1}{n - 1}.
\]

where $k_2^* = k_2 - k_1$ and $m^* = m - n k_1$.

\section{Derivations for Different Scenarios When $m < n$}
\label{section3}

When \(m < n\), distributing \(m\) identical balls into \(n\) distinct boxes inevitably results in some boxes remaining empty. Consequently, the constraint \(k_1 \leq k \leq k_2\) simplifies to \(0 \leq k \leq k_2\). The possible range of \(k_2\) can be categorized into three distinct cases: \(k_2 \geq m\), \(\lfloor \frac{m}{2} \rfloor \leq k_2 < m\), and \(\lceil \frac{m}{n} \rceil \leq k_2 < \lfloor \frac{m}{2} \rfloor\). Furthermore, since \(m < n\), the condition \(\lceil \frac{m}{n} \rceil \leq k_2 < \lfloor \frac{m}{2} \rfloor\) simplifies to \(1 \leq k_2 \leq \lfloor \frac{m}{2} \rfloor\).

\subsection{\(k_1 = 0\) and \(k_2 \geq m\)} 

As established in Section \ref{sec212}, when no upper bound is imposed on individual box capacity, the total number of valid distributions is given by: 
\[
\binom{m+n-1}{n-1}.
\]

\subsection{\(k_1 = 0\) and \(\lfloor \frac{m}{2} \rfloor \leq k_2 < m\)}

From the results in Section \ref{221}, when the maximum capacity constraint \(k_2\) satisfies \(\lfloor \frac{m}{2} \rfloor \leq k_2 < m\), the total number of valid distributions is:
\[
\binom{m+n-1}{n-1} - \binom{n}{1} \binom{m - k_2 + n - 2}{n - 1}.
\]

\subsection{\(k_1 = 0\) and \(1 \leq k_2 \leq \lfloor \frac{m}{2} \rfloor\)} 

Following the derivation in Section \ref{223}, when the upper bound on each box is limited to \(1 \leq k_2 \leq \lfloor \frac{m}{2} \rfloor\), the total number of distributions is:
\[
\binom{m+n-1}{n-1} - \sum_{\alpha=1}^{\lfloor \frac{m}{k_2+1} \rfloor}(-1)^{\alpha+1} \binom{n}{\alpha} \binom{m - \alpha (k_2+1) + n - 1}{n - 1}.
\]

\section{Conclusion}
\label{sec4}

This study extends classical combinatorial distribution theory by introducing and analyzing models with capacity constraints on urns. Through the development of a rigorous theoretical framework, we have systematically examined the effects of these constraints on the distribution of indistinguishable balls into distinct urns. The results highlight significant variations in distribution strategies, influenced by whether the constraints apply to the upper or lower bounds of urn capacities.

The analytical solutions derived in this work not only deepen our understanding of classical distribution problems but also introduce novel methodologies for addressing complex scenarios involving capacity constraints. These methodologies are particularly relevant to fields such as statistical physics and network theory, where precise resource allocation plays a critical role in optimizing systems and driving advancements.

Furthermore, this research establishes a foundation for future investigations into combinatorial problems with similar constraints. It encourages a reevaluation of existing models and opens new pathways for applying these refined theories to practical challenges in engineering, data science, and other disciplines requiring optimized allocation strategies.

In summary, this study expands the boundaries of combinatorial distribution theory by incorporating variable capacity constraints, addressing a fundamental mathematical challenge while enriching the toolkit available to researchers and practitioners tackling real-world problems across diverse fields.


\begin{thebibliography}{00}

\bibitem{feller1991introduction}
  W. Feller,
  An introduction to probability theory and its applications, Volume 2,
  John Wiley \& Sons, 81 (1991).

\bibitem{gazdzicki2020notes}
  M. Gazdzicki, M. I. Gorenstein, O. Savchuk, L. Tinti,
  Notes on statistical ensembles in the Cell Model,
  International Journal of Modern Physics E, 29 (2020) 2050060.
  https://doi.org/10.1142/S0218301320500603.
 

\bibitem{goodstein2007boolean}
  R. L. Goodstein,
  Boolean algebra,
  Courier Corporation, (2007).



\bibitem{gordon2002maxwell}
  B. L. Gordon, Maxwell–Boltzmann statistics and the metaphysics of modality, 
  \textit{Synthese} 133 (2002) 393–417. \newline https://doi.org/10.1023/A:1021360805193.


\bibitem{ko2008signature}
  K. Ko, T. G. Robertazzi,
  Signature search time evaluation in flat file databases,
  IEEE Transactions on Aerospace and Electronic Systems, 44 (2008) 493–502.
  https://doi.org/10.1109/TAES.2008.4560202.



\bibitem{kyong2012greedy}
  Y. Kyong, T. G. Robertazzi,
  Greedy signature processing with arbitrary location distributions: A divisible load framework,
  IEEE Transactions on Aerospace and Electronic Systems, 48 (2012) 3027–3041.
  https://doi.org/10.1109/TAES.2012.6324675.


\bibitem{ryser1963combinatorial}
  H. J. Ryser,
  Combinatorial mathematics,
  American Mathematical Society, 14 (1963).

\bibitem{ying2013signature}
  Z. Ying, T. G. Robertazzi,
  Signature searching in a networked collection of files,
  IEEE Transactions on Parallel and Distributed Systems, 25 (2013) 1339–1348.
  https://doi.org/10.1109/TPDS.2013.258.


\end{thebibliography}
\end{document}